\documentclass[a4paper,11pt]{article}
\usepackage{cite}
\usepackage{amsmath,amssymb,amsthm}
\usepackage{graphicx}
\usepackage{tikz}
\usepackage{booktabs}
\usepackage{longtable}
\usepackage{array}
\usepackage[ruled,vlined,linesnumbered]{algorithm2e}
\SetKwComment{Comment}{/* }{ */}
\usepackage{subcaption}
\usepackage{url}

\usepackage[colorlinks=true,linkcolor=blue,citecolor=blue,urlcolor=blue]{hyperref}
\newtheorem{theorem}{Theorem}[section]
\newtheorem{lemma}[theorem]{Lemma}
\newtheorem{corollary}[theorem]{Corollary}
\newtheorem{proposition}[theorem]{Proposition}
\newtheorem{definition}[theorem]{Definition}

\theoremstyle{remark}
\newtheorem{example}[theorem]{Example}
\newtheorem{remark}[theorem]{Remark}

\newcommand{\FF}{\mathbb{F}}

\DeclareMathOperator{\wt}{wt}

\DeclareMathOperator{\diag}{diag}

\title{Binary LCD Codes and Their Graph Representations}
\author{Keita~Ishizuka
\thanks{K. Ishizuka is with the Information Technology R\&D Center, Mitsubishi Electric Corporation, Kanagawa, Japan (e-mail: ishizuka.keita@ce.mitsubishielectric.co.jp).}
\thanks{This work was supported by JSPS KAKENHI Grant Number JP25K17290.}}
\date{}

\begin{document}
\maketitle

\begin{abstract}
We give a complete characterization of simple graphs whose adjacency matrices generate binary linear complementary dual (LCD) codes.
In particular, we completely characterize a distance-regular graph which yields an LCD code in terms of the intersection array parameters. This necessary and sufficient criterion strengthens the previously known sufficient conditions and unifies the cases of complete, Hamming, Johnson, and Grassmann graphs.
As further applications, we prove that non-isomorphic conference graphs with $q \equiv 1 \pmod 8$ yield inequivalent codes and we classify all simple graphs with idempotent adjacency matrices on at most $13$ vertices via mass formulas for binary LCD codes.
\end{abstract}

\medskip
\noindent\textbf{Keywords:} LCD codes, distance-regular graphs, strongly regular graphs.

\section{Introduction}

The code equivalence problem asks whether two linear codes are identical up to a coordinate permutation. This fundamental problem in coding theory has deep connections to computational complexity. Petrank and Roth \cite{petrank1997} established that code equivalence is at least as hard as graph isomorphism (GI), showing that if one can efficiently solve code equivalence, then one can also efficiently solve GI. This raised a natural question: does the converse hold? Is code equivalence computationally equivalent to graph isomorphism?
For LCD codes, Bardet et al.\ \cite{bardet2019} gave an affirmative answer: permutation code equivalence reduces to GI by exploiting a key structural property that the orthogonal projector $\Pi_C$ onto an LCD code is a symmetric matrix, enabling a polynomial-time reduction to GI.

The result of Bardet et al.\ suggests a deeper structural connection between LCD codes and graphs. Their reduction treats orthogonal projectors as general symmetric matrices, but in combinatorics, the most natural and well-studied objects are \emph{simple graphs}, i.e., graphs without loops, multiple edges, or directions. This motivates the investigation of which LCD codes correspond to simple graphs and what characterizes these codes. More broadly, such a correspondence could connect LCD codes to well-studied graph families and enable transfer of tools and techniques between graph theory and coding theory.

\subsection{Our Contributions}

This paper addresses these questions for binary LCD codes through a comprehensive framework connecting them with simple graphs via orthogonal projectors. Our main contributions are as follows.

\begin{enumerate}
\item \emph{A necessary and sufficient characterization of distance-regular graphs yielding LCD codes} (Theorem~\ref{thm:drg-main}). Given a distance-regular graph with intersection array $\{b_0, b_1, \ldots, b_{d-1}; c_1, c_2, \ldots, c_d\}$, the resulting binary code from its adjacency matrix is LCD if and only if the three parameters $b_0$, $a_1 = b_0 - b_1 - c_1$, and $c_2$ satisfy explicit parity conditions. This strengthens the previously known sufficient conditions of Key and Rodrigues~\cite{key2018} for strongly regular graphs to a complete characterization across all distance-regular graphs, and yields, as immediate corollaries, sharp criteria for the complete graphs $K_n$, the Hamming graphs $H(n,m)$, the Johnson graphs $J(n,k)$, the Grassmann graphs $J_q(n,k)$, and the cycle graphs $C_n$ to produce LCD codes (Theorem~\ref{thm:srg-from-drg} and Propositions~\ref{prop:cycle-graph}--\ref{prop:johnson}).

\item \emph{An equivalence-preserving bijection underlying the characterization} (Theorem~\ref{thm:main-binary}). To establish the criterion above and to enable systematic enumeration, we prove that binary even LCD codes correspond bijectively to simple graphs whose adjacency matrices are idempotent over $\FF_2$, with permutation equivalence of codes corresponding exactly to graph isomorphism. A combinatorial characterization in terms of vertex degrees and common neighbors is given in Corollary~\ref{cor:binary-char}.

\item \emph{Applications.} Combining the two results above, we prove that non-isomorphic conference graphs with $q \equiv 1 \pmod 8$ yield inequivalent binary LCD codes (Theorem~\ref{thm:conference}), providing the first theoretical explanation for observation~(i) of Haemers, Peeters, and van Rijckevorsel~\cite{haemers1999}. We further give a complete classification of all simple graphs with idempotent adjacency matrices on at most $13$ vertices, using mass formulas for binary LCD codes (Section~\ref{sec:computational}).
\end{enumerate}

\subsection{Related Work}
\emph{LCD codes.} LCD codes have found important applications in cryptography. Carlet and Guilley \cite{carlet2014complementary} demonstrated that binary LCD codes enable orthogonal direct sum masking (ODSM) as a countermeasure against side-channel attacks. A fundamental structural result by Carlet et al.~\cite{carlet2018} proved that any linear code over $\FF_q$ with $q > 3$ is equivalent to an LCD code; thus binary LCD codes form a proper subclass with distinct structural properties.
For classification purposes, Carlet et al.~\cite{carlet2019} developed mass formulas for LCD codes, which enable their classification without exhaustive computation. Araya and Harada \cite{araya2019} used these formulas to classify all inequivalent binary LCD codes of length $n \leq 13$, and Bouyuklieva, Bouykliev, and \"Ozbudak \cite{bouyuklieva2026} subsequently extended the enumeration up to length $n \leq 20$.
Given these practical and theoretical motivations, understanding the structure of binary LCD codes is important.

\emph{Code-graph connections.} Researchers have long studied constructions of codes from graph adjacency matrices. Brouwer and van Eijl \cite{brouwer1992} investigated the $p$-rank (the rank of a matrix over $\FF_p$) of adjacency matrices of strongly regular graphs (SRGs), establishing foundational connections between graph spectra and code dimensions. Haemers et al.\ \cite{haemers1999} constructed binary codes from SRGs and made notable computational observations: for $\text{srg}(25,12,5,6)$, every pair of non-isomorphic graphs yielded inequivalent codes, and the Paley graph achieved the largest minimum distance. They observed the same pattern for $\text{srg}(41,20,9,10)$ on their samples and conjectured it holds generally.
The work most closely related to ours is that of Key and Rodrigues~\cite{key2018}, who gave \emph{sufficient} conditions for an SRG to produce a binary LCD code in terms of specific relations among its parameters $(v,k,\lambda,\mu)$. Our Theorem~\ref{thm:drg-main} sharpens this in two directions: (i) it provides \emph{necessary and sufficient} conditions, expressed as three explicit parity conditions on the intersection array parameters $b_0$, $a_1$, $c_2$; and (ii) it applies to the strictly broader class of distance-regular graphs, with strongly regular graphs recovered as the diameter-$2$ case (Theorem~\ref{thm:srg-from-drg}).
The prior works above are one-directional, constructing codes from graphs and studying their properties. In contrast, our framework establishes an equivalence-preserving bijection that runs in both directions, enabling us to provide a theoretical explanation for observation (i) of Haemers et al.\ \cite{haemers1999} and to apply mass formulas \cite{carlet2019} for classifying graphs with idempotent adjacency matrices.

In summary, our work bridges LCD code theory and code-graph connections: graph-theoretic techniques such as distance-regularity yield new characterizations of LCD codes, while coding-theoretic tools such as mass formulas enable classification of graphs with idempotent adjacency matrices.

\subsection{Organization}

The remainder of this paper is organized as follows. Section~\ref{sec:prelim} establishes notation and reviews binary LCD codes, SRGs, and DRGs. Section~\ref{sec:correspondence} develops the bijection between binary even LCD codes and simple graphs with idempotent adjacency matrices. Section~\ref{sec:drg} characterizes DRGs via the three-term recurrence relation and then applies to well-known graph families. Section~\ref{sec:applications} presents applications to SRGs, provides a theoretical explanation for observation (i) of Haemers et al.\ on conference graphs, and provides a complete computational classification of LCD-derived graphs. Section~\ref{sec:conclusion} discusses implications and open problems.

\section{Preliminaries}\label{sec:prelim}

\subsection{Binary Linear Codes and LCD Codes}

Let $\FF_2$ denote the finite field of order $2$.
A binary $[n,k]$ code is a $k$-dimensional subspace of $\FF_2^n$. For a vector $x = (x_1,\ldots,x_n) \in \FF_2^n$, the \emph{support} is $\text{supp}(x) = \{i : x_i \neq 0\}$, the \emph{weight} is $\wt(x) = |\text{supp}(x)|$, and the \emph{minimum distance} of a code $C$ is $d(C) = \min\{\wt(c) : c \in C \setminus \{0\}\}$.
A code with parameters $[n,k,d]$ has length $n$, dimension $k$, and minimum distance $d$.
A binary $[n,k]$ code is \emph{optimal} if it attains the largest possible minimum distance among all binary $[n,k]$ codes.
A code is \emph{even} if all codewords have even weight.
A \emph{generator matrix} $G$ of an $[n,k]$ code $C$ is a $k \times n$ matrix whose rows form a basis for $C$.
Two codes $C_1$ and $C_2$ are \emph{(permutation) equivalent} if there exists a permutation matrix $P$ such that $C_2 = \{cP : c \in C_1\}$.

The \emph{Euclidean inner product} is $\langle x,y \rangle = \sum_{i=1}^n x_iy_i$, and the \emph{dual code} of $C$ is $C^\perp = \{x \in \FF_2^n : \langle x,c \rangle = 0 \text{ for all } c \in C\}$.
A code $C$ is \emph{linear complementary dual} (LCD) if $C \cap C^\perp = \{0\}$.
Equivalently, $C$ is LCD if and only if $\FF_2^n = C \oplus C^\perp$, since $\dim C + \dim C^\perp = n$ over $\FF_2$.
Throughout this paper, the central object associated with an LCD code is its orthogonal projector, defined as follows.

\begin{definition}[Orthogonal Projector of an LCD Code]\label{def:projector}
Let $C \subseteq \FF_2^n$ be an LCD code, so that $\FF_2^n = C \oplus C^\perp$. The \emph{orthogonal projector} of $C$ is the unique linear map $\Pi_C : \FF_2^n \to \FF_2^n$ satisfying $v\Pi_C = v$ for all $v \in C$ and $v\Pi_C = 0$ for all $v \in C^\perp$.
\end{definition}

Note that the LCD assumption $C \cap C^\perp = \{0\}$ is essential: it ensures that the two defining conditions of Definition~\ref{def:projector} are simultaneously consistent (if $v \in C \cap C^\perp$, then $v = v\Pi_C = 0$). The projector is represented by an $n \times n$ symmetric matrix over $\FF_2$ satisfying $\mathrm{Im}(\Pi_C) = C$ and $\Pi_C^2 = \Pi_C$. Given a generator matrix $G$ of $C$, Theorem~\ref{thm:massey} yields the explicit formula:

\begin{theorem}[Massey \cite{massey1992}]\label{thm:massey}
A binary code $C$ with generator matrix $G$ is LCD if and only if $GG^T$ is nonsingular over $\FF_2$. The orthogonal projector is given by
\[
\Pi_C = G^T(GG^T)^{-1}G.
\]
\end{theorem}

\subsection{Simple Graphs}

A \emph{simple graph} $\Gamma = (V,E)$ consists of a finite vertex set $V$ and edge set $E \subseteq \binom{V}{2}$ (no loops or multiple edges). The \emph{neighborhood} of a vertex $v$ is $N(v) = \{u \in V : \{u,v\} \in E\}$, and its \emph{degree} is $\deg(v) = |N(v)|$. A graph is \emph{$k$-regular} if every vertex has degree $k$.
The \emph{complement} of a graph $\Gamma=(V,E)$, denoted $\overline{\Gamma}$, is the graph on the same vertex set V with
$E(\overline{\Gamma})=\binom{V}{2}\setminus E.$
Equivalently, two distinct vertices are adjacent in $\overline{\Gamma}$ exactly when they are non-adjacent in $\Gamma$.
An \emph{adjacency matrix} $A$ of a graph on $n$ vertices is the $n \times n$ symmetric $(0,1)$-matrix with $A_{ij} = 1$ if and only if $\{i,j\} \in E$. Note that $A_{ii} = 0$ for all $i$ since there are no loops.
Two graphs are \emph{isomorphic} if their adjacency matrices $A$ and $A'$ satisfy $A' = P^TAP$ for some permutation matrix $P$.

A \emph{strongly regular graph} (SRG) with parameters $(v,k,\lambda,\mu)$, denoted $\text{srg}(v,k,\lambda,\mu)$, is a simple graph which has $v$ vertices, is $k$-regular, any two adjacent vertices have exactly $\lambda$ common neighbors, and any two non-adjacent vertices have exactly $\mu$ common neighbors.

\begin{lemma}[SRG Relations{~\cite[Theorem 9.1.2]{brouwerhaemers2012}}]\label{lem:srg-relations}
For an $\text{srg}(v,k,\lambda,\mu)$, the adjacency matrix $A$ satisfies
\[
A^2 = kI + \lambda A + \mu(J - I - A),
\]
where $I$ is the identity matrix and $J$ is the all-ones matrix.
\end{lemma}

Important families include (see~\cite[Chapter~9]{brouwerhaemers2012}):
\begin{itemize}
\item \emph{Complete graphs} $K_n$: The complete graph on $n$ vertices with all possible edges, having SRG parameters $(n,n-1,n-2,0)$.
\item \emph{Conference graphs}: A class of SRGs with parameters $(v,(v-1)/2,(v-5)/4,(v-1)/4)$ where $v \equiv 1 \pmod 4$. These graphs arise from conference matrices and include various construction methods.
\item \emph{Paley graphs}: A specific subfamily of conference graphs where $v = q$ is a prime power. Constructed using quadratic residues in $\FF_q$, the Paley graph is unique for each prime power $q \equiv 1 \pmod 4$.
\end{itemize}

A connected graph $\Gamma$ of diameter $d$ is called a \emph{distance-regular graph} (DRG) if there exist integers $b_i\ (0 \le i \le d-1)$ and $c_i$\ ($1 \leq i \leq d$) such that for any two vertices $x, y$ at distance $i$, the number of vertices at distance $i-1$ from $x$ and adjacent to $y$ equals $c_i$, and the number at distance $i+1$ from $x$ and adjacent to $y$ equals $b_i$.
The sequence $\{b_0, b_1, \ldots, b_{d-1}; c_1, c_2, \ldots, c_d\}$ is called the \emph{intersection array}.

\begin{lemma}[Three-Term Recurrence{~\cite[Proposition 2.71]{bannai2016}}]\label{lem:three-term}
Let $A_i$ denote the distance-$i$ matrix with $(A_i)_{xy} = 1$ if and only if $d(x,y) = i$. Setting $A=A_1$ and  $a_i = k - b_i - c_i$ where $k = b_0$ is the valency, we have the fundamental three-term recurrence:
\[
A \cdot A_i = b_{i-1}A_{i-1} + a_i A_i + c_{i+1}A_{i+1}.
\]
\end{lemma}

Important families include (see~\cite[Chapter~12]{brouwerhaemers2012}):
\begin{itemize}
\item \emph{Strongly regular graphs}: DRGs with diameter $2$, characterized by parameters $(v,k,\lambda,\mu)$ where $b_0 = k$, $a_1 = \lambda$, $c_2 = \mu$.
\item \emph{Cycle graphs} $C_n$:
We distinguish three cases: (i) For $n=2m\ (m \ge 1)$, $C_n$ has diameter $m$ with intersection array $b_0 = 2$, $b_i = 1\ (1\le i < m)$, $b_m = 0$ and $c_0 = 0$, $c_i = 1\ (1\le i < m)$, $c_m = 2$.
(ii) For $n=2m+1\ (m \ge 2)$, $C_n$ has diameter $m$ with intersection array $b_0 = 2$, $b_i = 1\ (1\le i < m)$, $b_m=0$ and $c_0=0$, $c_i = 1\ (1 \le i < m)$, $c_m = 1$.
(iii) If $n=3$, $C_n$ is the complete graph of order three $K_3$, which is an $\text{srg}(3,2,1,0)$.

\item \emph{Hamming graphs} $H(n,m)$: Vertices are $n$-tuples over an $m$-element alphabet, adjacent when differing in exactly one coordinate. The graph has diameter $n$ with intersection array $b_i = (m-1)(n-i)$ and $c_i = i$. For example, $H(3,2)$ is the $3$-dimensional cube graph.
\item \emph{Johnson graphs} $J(n,k)$: Vertices are $k$-subsets of an $n$-set, adjacent when they share $k-1$ elements. The graph has diameter $\min(k, n-k)$ with $b_i = (k-i)(n-k-i)$ and $c_i = i^2$. The complement of $J(5,2)$ is the famous Petersen graph.
\item \emph{Grassmann graphs} $J_q(n,k)$: 
Vertices are $k$-dimensional subspaces of $\FF_q^n$, adjacent if their intersection has dimension $k-1$.
The graph has diameter $\min(k,n-k)$ with intersection arrays
$b_i = q^{2i+1}\binom{k-i}{1}_q \binom{n-k-i}{1}_q$ and $c_i = \binom{i}{1}_q^2$.
Here the Gaussian binomial coefficient is defined by
$\binom{n}{k}_q = \frac{(q^n-1)(q^{n-1}-1)\cdots(q^{n-k+1}-1)}
{(q^k-1)(q^{k-1}-1)\cdots(q-1)}.$
\end{itemize}

\section{LCD-Graph Correspondence}\label{sec:correspondence}

We establish the fundamental bijection between binary even LCD codes and simple graphs with idempotent adjacency matrices.
We begin with the following folklore lemma, and include a proof for completeness.

\begin{lemma}[Characterization of Binary Even Code]\label{lem:evencode}
A binary code is even if and only if every row of a generator matrix has even weight.
\end{lemma}
\begin{proof}
Let $x,y\in \FF_2^n$. It is well known~\cite[Theorem 1.4.3]{huffman2010} that
\[
\wt(x+y)=\wt(x)+\wt(y)-2\wt(x\cap y),
\]
where $x\cap y$ denotes the vector having $1$s precisely in the coordinates where both $x$ and $y$ have $1$s. In particular,
\[
\wt(x+y)\equiv \wt(x)+\wt(y)\pmod 2.
\]
Hence the parity of the weight is preserved under addition: the sum of two even vectors is even.

Let $G$ be a generator matrix whose rows all have even weight. Every codeword is a linear combination of the rows of $G$, i.e., a sum of some subset of these rows. By repeatedly applying the above parity relation, any such sum has even weight. Therefore the code is even.

Conversely, if the code is even, then every codeword has even weight. Since each row of a generator matrix is itself a codeword, every row must have even weight.
\end{proof}

\begin{theorem}[Binary LCD-Graph Correspondence]\label{thm:main-binary}
There is a bijection between:
\begin{itemize}
\item Binary even LCD codes of length $n$;
\item Simple graphs on $n$ vertices whose adjacency matrix $A$ is idempotent over $\FF_2$.
\end{itemize}
This bijection preserves equivalence: two codes are permutation equivalent if and only if the corresponding graphs are isomorphic.
\end{theorem}

\begin{proof}
Let $A = \Pi_C$ be the orthogonal projector, and let $r_i$ be the $i$-th row of $A$.
We first claim that $A_{ii} \equiv \wt(r_i) \pmod 2$.
Indeed, for binary codes, the Euclidean inner product satisfies $\langle v, v \rangle = \sum_i v_i^2 = \sum_i v_i \equiv \wt(v) \pmod 2$.
Since $A$ is idempotent, $(A^2)_{ii}=A_{ii}$. But $(A^2)_{ii} = \langle r_i, r_i \rangle \equiv \wt(r_i) \pmod 2$. Therefore $A_{ii} \equiv \wt(r_i) \pmod 2$, as required.

\emph{From codes to graphs:} Let $C$ be a binary even LCD code with orthogonal projector $A$. Since every codeword in $C$ has even weight and the rows of $A$ generate $C$, each row $r_i$ satisfies $\wt(r_i) \equiv 0 \pmod 2$. Thus $A_{ii} = 0$ for all $i$. Since $A$ is symmetric (being an orthogonal projector) and has zero diagonal, it is the adjacency matrix of a simple graph. The idempotence $A^2 = A$ holds over $\FF_2$ by definition of the projector.

\emph{From graphs to codes:} Let $\Gamma$ be a simple graph whose adjacency matrix $A$ is idempotent over $\FF_2$. Since $A$ is symmetric and idempotent, it is an orthogonal projector. By Theorem~\ref{thm:massey}, the row span $C$ of $A$ is an LCD code. Since $A_{ii} = 0$ (no loops), we have $\wt(r_i) \equiv 0 \pmod 2$ for each row $r_i$.
By Lemma~\ref{lem:evencode}, $C$ is an even code.

\emph{Equivalence preservation:} If two LCD codes $C$ and $C'$ are permutation equivalent, their generator matrices satisfy $G' = MG P$ for some invertible matrix $M$ and permutation matrix $P$. For the orthogonal projectors:
\begin{align*}
A'
&= (G')^T((G')(G')^T)^{-1}G'\\
&= P^T G^T M^T((MGP)(P^T G^T M^T))^{-1}(MGP)\\
&= P^T G^T M^T (M^T)^{-1} (GG^T)^{-1} (M)^{-1} MGP\\
&= P^TAP.
\end{align*}
In the third equality, we use the fact that $GG^T$ is invertible since $C$ is LCD (Theorem~\ref{thm:massey}).
Therefore, the corresponding graphs are isomorphic.
Conversely, if $A' = P^TAP$ for some permutation matrix $P$, then the codes are permutation equivalent since their projectors determine the codes as their row spans.
\end{proof}

The graphs in Theorem~\ref{thm:main-binary} have a nice combinatorial characterization:

\begin{corollary}[Binary Combinatorial Characterization]\label{cor:binary-char}
A simple graph $\Gamma$ yields a binary even LCD code if and only if:
\begin{itemize}
\item Every vertex has even degree;
\item Any two adjacent vertices have an odd number of common neighbors;
\item Any two non-adjacent vertices have an even number of common neighbors.
\end{itemize}
\end{corollary}

\begin{proof}
Let $A$ be the adjacency matrix. The idempotence condition over $\FF_2$ means $(A^2)_{ij} = A_{ij}$ for all $i,j$. For $i = j$, $(A^2)_{ii}$ counts the degree of vertex $i$, while $A_{ii} = 0$. Thus the degree must be even. For $i \neq j$, $(A^2)_{ij}$ counts walks of length $2$ from vertex $i$ to vertex $j$, which equals the number of common neighbors. The result follows.
\end{proof}

\begin{remark}\label{rem:excluded-families}
The conditions in Corollary~\ref{cor:binary-char} immediately exclude many common graph families:
\begin{itemize}
\item No tree can yield an LCD code (trees have leaves with degree $1$);
\item No bipartite graph can yield an LCD code (adjacent vertices in bipartite graphs have no common neighbors);
\item More generally, any graph yielding an LCD code must contain triangles (since adjacent vertices need an odd number of common neighbors).
\end{itemize}
\end{remark}

\section{Distance-Regular Graph Characterization}\label{sec:drg}

We now present our main theoretical framework that unifies the study of LCD codes from graph families. Distance-regular graphs provide the right level of generality—more structured than arbitrary graphs, yet broad enough to encompass important families. For a DRG with diameter $d$, we have distance matrices $A_0 = I, A_1 = A, A_2, \ldots, A_d$ where $(A_i)_{xy} = 1$ if and only if the distance between vertices $x$ and $y$ equals $i$. The key property is the three-term recurrence relation (Lemma~\ref{lem:three-term}):
\[
A \cdot A_i = b_{i-1}A_{i-1} + a_i A_i + c_{i+1}A_{i+1}.
\]
In particular, for $i = 1$, we have
\[
A^2 = b_0 I + a_1 A + c_2 A_2,
\]
where $b_0 = k$ is the valency and $a_1 = k - b_1 - c_1$.

\begin{theorem}[DRG Characterization]\label{thm:drg-main}
A distance-regular graph $\Gamma$ with intersection array $\{b_0, b_1, \ldots, b_{d-1}; c_1, c_2, \ldots, c_d\}$ yields a binary even LCD code if and only if (i) $b_0 \equiv 0 \pmod 2$; (ii) $a_1 \equiv 1 \pmod 2$ where $a_1 = k - b_1 - c_1$; and (iii) $c_2 \equiv 0 \pmod 2$.
\end{theorem}

\begin{proof}
From the three-term recurrence relation, $A^2 = b_0 I + a_1 A + c_2 A_2$.
It is easy to see that $A^2 \equiv A \pmod 2$ holds if the coefficient of $I$ is $0 \pmod 2$ (i.e., $b_0 \equiv 0 \pmod 2$), the coefficient of $A$ is $1 \pmod 2$ (i.e., $a_1 \equiv 1 \pmod 2$), and the coefficient of $A_2$ is $0 \pmod 2$ (i.e., $c_2 \equiv 0 \pmod 2$).

Conversely, if any of these conditions fails, then $A^2 \not\equiv A \pmod 2$. The key observation is that $I$, $A$, and $A_2$ are (0,1)-matrices with disjoint support: $I$ has 1s only on the diagonal, $A$ has 1s only at positions $(i,j)$ where $d(i,j) = 1$, and $A_2$ has 1s only at positions $(i,j)$ where $d(i,j) = 2$. Since these matrices have 1s in completely different positions, they cannot cancel each other over $\FF_2$. If $b_0 \equiv 1 \pmod 2$, then $A^2$ has 1s on the diagonal (from $I$), but $A$ has zero diagonal, so $A^2 \neq A$. If $a_1 \equiv 0 \pmod 2$, then the $A$ component vanishes in $A^2 \pmod 2$, but we need exactly $A$ to remain. If $c_2 \equiv 1 \pmod 2$, then $A^2$ has 1s at distance-2 positions (from $A_2$), which $A$ doesn't have, so $A^2 \neq A$. Therefore, these conditions are both necessary and sufficient.
\end{proof}

\subsection{Applications to Special Graph Classes}

We now show how our main DRG theorem specializes to important graph families.

\begin{theorem}[SRG Characterization]\label{thm:srg-from-drg}
A strongly regular graph with parameters $(v, k, \lambda, \mu)$ yields a binary even LCD code if and only if (i) $k \equiv 0 \pmod 2$; (ii) $\lambda \equiv 1 \pmod 2$; and (iii) $\mu \equiv 0 \pmod 2$.
\end{theorem}

\begin{proof}
For an SRG, the DRG parameters specialize to $b_0 = k$, $a_1 = \lambda$, and $c_2 = \mu$. The result follows directly from Theorem~\ref{thm:drg-main}.
\end{proof}

\begin{example}[Paley Graphs]\label{ex:paley}
The Paley graph $P(9)$ is a strongly-regular graph with parameters $(9,4,1,2)$ and satisfies the hypothesis of Theorem~\ref{thm:srg-from-drg}.
The binary code generated by the rows of its adjacency matrix is a binary $[9,4,4]$ even LCD code.
Its generator matrix is
\[
G = \left(\begin{array}{rrrrrrrrr}
1 & 0 & 0 & 0 & 1 & 0 & 1 & 1 & 0 \\
0 & 1 & 0 & 0 & 0 & 1 & 0 & 1 & 1 \\
0 & 0 & 1 & 0 & 1 & 1 & 0 & 1 & 0 \\
0 & 0 & 0 & 1 & 1 & 1 & 1 & 1 & 1
\end{array}\right).
\]
One checks that the support of the orthogonal projector $A=G^T(GG^T)^{-1}G$ coincides with an adjacency matrix of $P(9)$.
According to Grassl's tables~\cite{grassl}, this code is optimal.
A detailed study of the Paley graph is given in the next section.
\end{example}

\begin{corollary}[Complete Graph Characterization]\label{cor:complete-graph}
The complete graph $K_n$ yields a binary even LCD code if and only if $n$ is odd.
\end{corollary}

\begin{proof}
$K_n$ is strongly regular with parameters $(n, n-1, n-2, 0)$. Applying Theorem~\ref{thm:srg-from-drg}, we need $k = n-1 \equiv 0 \pmod 2$, $\lambda = n-2 \equiv 1 \pmod 2$, and $\mu = 0 \equiv 0 \pmod 2$. The first two conditions both require $n$ odd, while the third is always satisfied. Thus $K_n$ yields an LCD code if and only if $n$ is odd.
\end{proof}

\begin{proposition}[Cycle Graph Characterization]\label{prop:cycle-graph}
Among all cycle graphs, only $C_3 = K_3$ yields a binary even LCD code.
\end{proposition}

\begin{proof}
For $n \geq 4$, the cycle graph $C_n$ is a distance-regular graph with parameters $b_0 = 2$ and $a_1 = 0$. By Theorem~\ref{thm:drg-main}, the condition $a_1 \equiv 1 \pmod 2$ is necessary, but $a_1 = 0 \equiv 0 \pmod 2$ fails this requirement. Thus $C_n$ does not yield an LCD code for $n \geq 4$. The case $C_3 = K_3$ yields an LCD code by Corollary~\ref{cor:complete-graph}.
\end{proof}

\begin{proposition}[Hamming Graph Characterization]
The Hamming graph $H(n,m)$ yields a binary even LCD code if and only if $m$ is odd.
\end{proposition}

\begin{proof}
The Hamming graph $H(n,m)$ has intersection array with
$$
b_i = (n-i)(m-1),\ c_i = i,\ a_i=(m-2)i.
$$
From these equations we verify that $b_0 = (m-1)n$, $a_1 = m-2$, $c_2 = 2$.
By Theorem~\ref{thm:drg-main}, $H(n,m)$ yields an LCD code if and only if $b_0 = (m-1)n \equiv 0 \pmod 2$, $a_1 = (m-2) \equiv 1 \pmod 2$, and $c_2 = 2 \equiv 0 \pmod 2$ (always satisfied). This requires $(m-1)n$ to be even and $m-2$ to be odd, which holds when $m$ is odd.
\end{proof}

\begin{proposition}[Johnson Graph Characterization]\label{prop:johnson}
The Johnson graph $J(n,k)$ yields a binary even LCD code if and only if $n$ is odd.
\end{proposition}

\begin{proof}
The Johnson graph $J(n,k)$ has intersection array with
$$
b_i = (k-i)(n-k-i),\ c_i=i^2,\ a_i=i(n-2i).
$$
From these equations we verify that $b_0 = k(n-k)$, $a_1 = n-2$, $c_2 = 4$.
By Theorem~\ref{thm:drg-main}, $J(n,k)$ yields an LCD code if and only if $b_0 = k(n-k) \equiv 0 \pmod 2$, $a_1 = n-2 \equiv 1 \pmod 2$, and $c_2 = 4 \equiv 0 \pmod 2$ (always satisfied). This requires $n$ odd.
\end{proof}
\begin{example}[The Petersen Graph]\label{ex:petersen}
The complement of the Petersen graph is isomorphic to $J(5,2)$ and it satisfies the hypothesis of Proposition~\ref{prop:johnson}.
The binary code generated by the rows of its adjacency matrix is a binary $[10,4,4]$ even LCD code.
A generator matrix is
\[
G = \left(\begin{array}{rrrrrrrrrr}
1 & 0 & 0 & 0 & 1 & 1 & 1 & 0 & 1 & 1 \\
0 & 1 & 0 & 0 & 1 & 0 & 0 & 1 & 1 & 0 \\
0 & 0 & 1 & 0 & 1 & 1 & 1 & 0 & 0 & 0 \\
0 & 0 & 0 & 1 & 1 & 1 & 0 & 1 & 1 & 1
\end{array}\right).
\]
According to Grassl's tables~\cite{grassl}, this code is optimal.
\end{example}

\begin{proposition}[Grassmann Graph Characterization]
The Grassmann graph $J_q(n,k)$ with $q$ odd yields a binary even LCD code if and only if $n$ is odd. If $q$ is even, then $J_q(n,k)$ never yields a binary even LCD code.
\end{proposition}
\begin{proof}
The Grassmann graph $J_q(n,k)$ has intersection array with
$$
b_i = q^{2i+1}\binom{k-i}{1}_q \binom{n-k-i}{1}_q,\ c_i = \binom{i}{1}_q^2.
$$
A straightforward calculation shows that
\begin{align*}
b_i
&= q^{2i+1} \binom{k-i}{1}_q \binom{n-k-i}{1}_q\\
&= \frac{q^{2i+1}(q^{k-i}-1)(q^{n-k-i}-1)}{(q-1)(q-1)}\\
&= q^{2i+1} (q^{k-i-1}+q^{k-i-2}+\dots+1)(q^{n-k-i-1}+q^{n-k-i-2}+\dots+1)\\
&\equiv (k-i)(n-k-i) \pmod 2.
\end{align*}
In the fourth equality, we used that $q \equiv 1 \pmod 2$.
In a similar manner, we see that
$$
c_i
= \binom{i}{1}_q^2
= (\frac{q^i-1}{q-1})^2
= (q^{i-1}+ \dots + q+1)^2
\equiv i \pmod 2.
$$
Using these equations, we obtain
\begin{align*}   
a_1
&= k - b_1 - c_1\\
&\equiv k(n-k) - (k-1)(n-k-1) - 1 \pmod 2\\
&\equiv kn - k^2 - (kn - k^2 - k - n + k + 1) - 1 \pmod 2\\
&\equiv n \pmod 2.
\end{align*}
By Theorem~\ref{thm:drg-main}, $J_q(n,k)$ yields an LCD code if and only if $b_0 \equiv 0 \pmod 2$, $a_1 \equiv 1 \pmod 2$, and $c_2 \equiv 0 \pmod 2$. With the above reductions, these become $k(n-k) \equiv 0$, $n \equiv 1$, and $c_2 \equiv 0 \pmod 2$, respectively. Since $c_2 \equiv 0$ is automatically satisfied, the conditions reduce to $n$ being odd.

If $q$ is even, the above reductions do not apply.
Instead, since $q \equiv 0 \pmod 2$, we have $q^j \equiv 0 \pmod 2$ for $j \geq 1$, so $\binom{r}{1}_q = q^{r-1}+\cdots+q+1 \equiv 1 \pmod 2$ for all $r \geq 1$. Therefore $c_2 = \binom{2}{1}_q^2 = (q+1)^2 \equiv 1 \pmod 2$, violating condition (iii) of Theorem~\ref{thm:drg-main}. Thus $J_q(n,k)$ never yields a binary even LCD code when $q$ is even.
\end{proof}

\section{Applications and Computational Results}\label{sec:applications}

This section demonstrates applications of our main results. In Section~\ref{sec:srg-applications}, we apply Theorem~\ref{thm:main-binary} and Theorem~\ref{thm:srg-from-drg} to conference graphs, providing a theoretical explanation for observations of Haemers et al.\ \cite{haemers1999}. In Section~\ref{sec:computational}, we classify simple graphs with idempotent adjacency matrices by leveraging mass formulas for LCD codes.

\subsection{Applications to Strongly Regular Graphs}\label{sec:srg-applications}

In their seminal computational study, Haemers, Peeters, and Rijckevorsel \cite{haemers1999} investigated binary codes from strongly regular graphs and made notable observations. For $\text{srg}(25,12,5,6)$, which had been completely classified, they verified two properties:
\begin{itemize}
\item[(i)] Every pair of non-isomorphic graphs yielded a pair of inequivalent codes;
\item[(ii)] The Paley graph achieved the largest minimum distance among all codes.
\end{itemize}
For $\text{srg}(41,20,9,10)$, whose complete classification remains open, they observed the same properties on their computational samples.

In this subsection, we apply our main results to provide a proof of observation (i) for all conference graphs with $q \equiv 1 \pmod 8$, and computationally investigate observation (ii) in the case $q = 41$.

\subsubsection{Conference Graphs}

Conference graphs are strongly regular graphs with parameters $(q,(q-1)/2,(q-5)/4,(q-1)/4)$ where $q$ is a prime power with $q \equiv 1 \pmod 4$. The Paley graph is the primary example.
Our theoretical framework provides new insights into the relationship between conference graphs and LCD codes:

\begin{lemma}[$p$-Rank of Conference Graphs{\cite[Section 4 (iii)]{brouwer1992}}]\label{lem:p-rank}
Let $\Gamma$ be a conference graph with $(v,k,\lambda,\mu)=(q,(q-1)/2,(q-5)/4,(q-1)/4)$ with adjacency matrix $A$.
If $\mu \equiv 0 \pmod p$, then its $p$-rank (i.e., the rank of $A$ over $\mathbb{F}_p$) equals $(q-1)/2$.
\end{lemma}

\begin{theorem}[Conference Graph Correspondence]\label{thm:conference}
A conference graph $(q,(q-1)/2,(q-5)/4,(q-1)/4)$ yields a binary even LCD code if and only if $q \equiv 1 \pmod 8$.
Furthermore, these codes have parameters $[q, (q-1)/2]$, and two codes are permutation equivalent if and only if the corresponding conference graphs are isomorphic.
\end{theorem}

\begin{proof}
We first note that the dimensions of the resulting codes equal the $2$-ranks of their adjacency matrices, which are $(q-1)/2$ by Lemma~\ref{lem:p-rank}.
Thus, all codes have the same lengths and dimensions.
By Lemma~\ref{lem:srg-relations}, the adjacency matrix $A$ satisfies
\[
A^2 = \frac{q-1}{2}I + \frac{q-5}{4}A + \frac{q-1}{4}(J-I-A).
\]
By Theorem~\ref{thm:srg-from-drg}, these graphs yield binary even LCD codes if and only if $(q-1)/2 \equiv 0 \pmod 2$, $(q-5)/4 \equiv 1 \pmod 2$, and $(q-1)/4 \equiv 0 \pmod 2$, and the condition is equivalent to $q \equiv 1 \pmod 8$.
The last assertion is an immediate consequence of Theorem~\ref{thm:main-binary}.
\end{proof}

Note that both $\text{srg}(25,12,5,6)$ and $\text{srg}(41,20,9,10)$ are conference graphs with $q = 25$ and $q = 41$, respectively, satisfying $q \equiv 1 \pmod 8$.
Theorem~\ref{thm:conference} provides a theoretical explanation for observation~(i) of Haemers et al.: non-isomorphic conference graphs yield inequivalent codes by the equivalence-preserving bijection of Theorem~\ref{thm:main-binary}.

\subsubsection{Optimality of Paley Graphs}
We now address observation (ii).
By Lemma~\ref{lem:p-rank}, all conference graphs with the same parameter $q \equiv 1 \pmod 8$ yield codes of the same length $q$ and dimension $(q-1)/2$. Thus the natural question is which graph achieves the largest minimum distance.

We examined all $7{,}152$ known non-isomorphic graphs with parameters $\text{srg}(41,20,9,10)$ listed in \cite{maksimovic2022}.
Among these graphs, only the Paley graph achieves minimum distance $10$, while all others yield codes with $4 \leq d \leq 9$.
By Grassl's tables~\cite{grassl}, the largest minimum distance of $[41,20]$ codes is $10$.
Therefore, observation (ii) holds for $\text{srg}(41,20,9,10)$.
Note that this does not rule out the existence of another conference graph yielding a code with minimum distance $10$, since a complete classification of $\text{srg}(41,20,9,10)$ has not been completed yet.

A general conjecture --- that among all conference graphs with parameters $\text{srg}(q,(q-1)/2,(q-5)/4,(q-1)/4)$ where $q \equiv 1 \pmod 8$, the Paley graph always yields the code with largest minimum distance --- remains open for $q > 41$.

\subsection{Classification of LCD-Derived Graphs of Small Orders}\label{sec:computational}

Carlet et al.~\cite{carlet2019} developed mass formulas for LCD codes, which enable systematic classification without exhaustive search over all generator matrices. No analogous mass formula exists for graphs with idempotent adjacency matrices. Our bijection (Theorem~\ref{thm:main-binary}) bridges this gap: by classifying LCD codes via mass formulas and then applying the correspondence, we can classify the associated graphs without direct enumeration over all graphs.
We demonstrate this approach using the database of Araya and Harada \cite{araya2019}, who classified binary LCD codes of length $n \leq 13$.
We note that Bouyuklieva, Bouykliev, and \"Ozbudak \cite{bouyuklieva2026} subsequently extended the enumeration of inequivalent binary LCD codes (with dual distance at least~$2$) up to length $n \leq 20$. Their results, however, report only the \emph{number} of inequivalent codes for each $(n,k)$; the underlying code database is not made publicly available. Since computing graph-theoretic invariants such as edge count, diameter, automorphism group order, and family membership requires explicit access to each code's generator matrix, we restrict our classification to the range covered by the database of \cite{araya2019}.

\subsubsection{Classification Methods Using Mass Formula}

The mass formula for binary LCD codes, established by Carlet et al.~\cite{carlet2019}, states that
\[
T_2(n,k) = \sum_{C \in \mathcal{B}_{n,k}} \frac{n!}{|\mathrm{Aut}(C)|},
\]
where $\mathcal{B}_{n,k}$ denotes the set of all inequivalent binary LCD $[n,k]$ codes and $\mathrm{Aut}(C)$ is the automorphism group of $C$. The value $T_2(n,k)$ can be computed theoretically using Gaussian binomial coefficients without explicitly finding $\mathcal{B}_{n,k}$, which can be found in~\cite{carlet2019}.
Araya and Harada \cite{araya2019} used this formula to classify binary LCD codes of lengths up to $13$. Their method proceeds as follows: enumerate candidate codes with generator matrices of the form $(I_k \mid A)$, test equivalence to build a set $\mathcal{B}'_{n,k}$ of inequivalent codes, and verify completeness by checking that $\sum_{C \in \mathcal{B}'_{n,k}} \frac{n!}{|\mathrm{Aut}(C)|} = T_2(n,k)$.

Our bijection (Theorem~\ref{thm:main-binary}) enables application of this classification to graphs: each binary even LCD code corresponds to a unique simple graph with idempotent adjacency matrix. Since the bijection preserves equivalence, inequivalent codes yield non-isomorphic graphs, and the mass formula verification ensures completeness of the graph classification.

\subsubsection{Computational Methods and Algorithms}

Our computational approach systematically processes the Araya-Harada database \cite{araya2019} of LCD codes. The database contains LCD codes with dimension $k \leq \lfloor n/2 \rfloor$, utilizing the fact that if $C$ is LCD, then so is its dual $C^\perp$. To obtain a complete classification of graphs from all LCD codes (including those with $k > n/2$), we examine both each code and its dual.

For each code $C$ in the database, we:
\begin{enumerate}
\item Compute the orthogonal projector $A = G^T(GG^T)^{-1}G$ where $G$ is a generator matrix of $C$;
\item Check the diagonal of $A$:
   \begin{itemize}
   \item If $\diag(A) = 0$ (all-zero): The code $C$ yields graph $\Gamma(A)$ directly;
   \item If $\diag(A) = \mathbf{1}$ (all-one): The dual code $C^\perp$ (with dimension $n-k \geq n/2$) yields graph $\Gamma(I+A)$;
   \item If diagonal is mixed: No graph correspondence exists.
   \end{itemize}
\item For valid cases, construct and analyze the resulting graph;
\item Store the pair $(C, \Gamma)$ with all computed properties.
\end{enumerate}

This dual examination ensures complete coverage: codes with $k \leq n/2$ are checked directly for all-zero diagonal, while codes with $k > n/2$ are captured through their duals when they have all-one diagonal.
For each graph, we compute the number of edges, diameter, number of connected components, and automorphism group order. All computations were performed using SageMath~\cite{sagemath}.

\subsubsection{Simple Graphs from Binary LCD Codes}\label{sec:binary-classification}

From 22,213 binary LCD codes with $n \leq 13$, we obtained 1,208 non-isomorphic graphs across 58 parameter sets $[n,k,d]$.
Table~\ref{tab:binary-statistics} summarizes the classification, where $|E|$ denotes edge count, ``Diam'' diameter, ``C'' components, and $|$Aut$|$ automorphism group order.
We use ``K'' for thousands and ``M'' for millions. Parameters containing well-known graph families are marked with asterisks~(*).

\begin{table}[h]
\centering
\caption{Statistics of graphs from binary LCD codes}
\label{tab:binary-statistics}
\scriptsize
\setlength{\tabcolsep}{3pt}
\begin{tabular}{cccccc|cccccc}
\hline
$[n,k,d]$ & \# & $|E|$ & Diam & C & $|$Aut$|$ & $[n,k,d]$ & \# & $|E|$ & Diam & C & $|$Aut$|$ \\
\hline
$[3,2,2]$* & 1 & 3 & 1 & 1 & 6 & $[11,2,4]$ & 3 & 15-31 & 2 & 1-5 & 1K-30K \\
$[4,2,2]$ & 1 & 3 & 1 & 2 & 6 & $[11,2,6]$ & 3 & 27-39 & 2 & 1-3 & 2K-28K \\
$[5,2,2]$ & 2 & 3-7 & 1-2 & 1-3 & 12 & $[11,4,2]$ & 39 & 6-46 & 1-4 & 1-7 & 16-120K \\
$[5,4,2]$* & 1 & 10 & 1 & 1 & 120 & $[11,4,4]$ & 8 & 18-34 & 2-3 & 1-3 & 16-144 \\
$[6,2,2]$ & 2 & 3-7 & 1-2 & 2-4 & 12-36 & $[11,6,2]$ & 46 & 9-49 & 1-4 & 1-5 & 16-120K \\
$[6,4,2]$ & 2 & 6-10 & 1 & 2 & 72-120 & $[11,6,4]$ & 1 & 25 & 2 & 1 & 120 \\
$[7,2,2]$ & 3 & 3-11 & 1-2 & 1-5 & 24-240 & $[11,8,2]$ & 11 & 16-52 & 1-2 & 1-3 & 1K-725K \\
$[7,2,4]$ & 1 & 15 & 2 & 1 & 72 & $[11,10,2]$* & 1 & 55 & 1 & 1 & 39M \\
$[7,4,2]$ & 4 & 6-18 & 1-2 & 1-3 & 24-240 & $[12,2,2]$ & 5 & 3-19 & 1-2 & 2-10 & 28K-2M \\
$[7,6,2]$* & 1 & 21 & 1 & 1 & 5K & $[12,2,4]$ & 3 & 15-31 & 2 & 2-6 & 4K-30K \\
$[8,2,2]$ & 3 & 3-11 & 1-2 & 2-6 & 72-720 & $[12,2,6]$ & 3 & 27-39 & 2 & 2-4 & 7K-28K \\
$[8,2,4]$ & 1 & 15 & 2 & 2 & 72 & $[12,4,2]$ & 58 & 6-46 & 1-4 & 1-8 & 32-604K \\
$[8,4,2]$ & 6 & 6-18 & 1-4 & 1-4 & 16-720 & $[12,4,4]$ & 18 & 18-42 & 2-4 & 1-4 & 16-432 \\
$[8,6,2]$ & 2 & 13-21 & 1 & 2 & 720-5K & $[12,6,2]$ & 113 & 9-49 & 1-5 & 1-6 & 4-604K \\
$[9,2,2]$ & 4 & 3-15 & 1-2 & 1-7 & 288-10K & $[12,6,4]$ & 4 & 25-33 & 2-3 & 1-2 & 8-120 \\
$[9,2,4]$ & 2 & 15-23 & 2 & 1-3 & 144-720 & $[12,8,2]$ & 33 & 12-52 & 1-4 & 1-4 & 96-2M \\
$[9,2,6]$* & 1 & 27 & 2 & 1 & 1K & $[12,10,2]$ & 3 & 31-55 & 1 & 2 & 604K-39M \\
$[9,4,2]$ & 13 & 6-30 & 1-4 & 1-5 & 16-2K & $[13,2,2]$ & 6 & 3-23 & 1-2 & 1-11 & 172K-79M \\
$[9,4,4]$* & 1 & 18 & 2 & 1 & 72 & $[13,2,4]$ & 4 & 15-39 & 2 & 1-7 & 17K-2M \\
$[9,6,2]$ & 7 & 9-33 & 1-2 & 1-3 & 144-10K & $[13,2,6]$ & 5 & 27-51 & 2 & 1-5 & 17K-604K \\
$[9,8,2]$* & 1 & 36 & 1 & 1 & 362K & $[13,2,8]$ & 1 & 55 & 2 & 1 & 172K \\
$[10,2,2]$ & 4 & 3-15 & 1-2 & 2-8 & 1K-30K & $[13,4,2]$ & 103 & 6-62 & 1-4 & 1-9 & 48-8M \\
$[10,2,4]$ & 2 & 15-23 & 2 & 2-4 & 432-720 & $[13,4,4]$ & 48 & 18-66 & 2-4 & 1-5 & 16-31K \\
$[10,2,6]$ & 1 & 27 & 2 & 2 & 1K & $[13,4,6]$ & 2 & 46-50 & 2 & 1 & 48-384 \\
$[10,4,2]$ & 20 & 6-30 & 1-4 & 1-6 & 16-14K & $[13,6,2]$ & 369 & 9-69 & 1-6 & 1-7 & 4-3M \\
$[10,4,4]$* & 3 & 18-30 & 2-3 & 1-2 & 16-120 & $[13,6,4]$ & 37 & 25-57 & 2-3 & 1-3 & 2-720 \\
$[10,6,2]$ & 14 & 9-33 & 1-4 & 1-4 & 48-30K & $[13,8,2]$ & 153 & 12-72 & 1-4 & 1-5 & 16-8M \\
$[10,8,2]$ & 3 & 20-36 & 1 & 2 & 28K-362K & $[13,10,2]$ & 16 & 23-75 & 1-2 & 1-3 & 17K-79M \\
$[11,2,2]$ & 5 & 3-19 & 1-2 & 1-9 & 5K-725K & $[13,12,2]$* & 1 & 78 & 1 & 1 & 6227M \\
\hline
\end{tabular}
\end{table}

Among the 1,208 graphs, we identify 44 that belong to well-known graph families: 6 complete graphs, 36 complete multipartite graphs, and 2 strongly regular graphs. Table~\ref{tab:known-families} lists these known families. Notably, the complete graphs $K_n$ for odd $n \leq 13$ all appear, consistent with Corollary~\ref{cor:complete-graph}. The complete multipartite graphs $K_{n_1,\ldots,n_r}$ arise when the partition sizes sum to $n$ with appropriate parity conditions.

\begin{table}[h]
\centering
\caption{Known graph families from binary LCD codes ($n \leq 13$)}
\label{tab:known-families}
\footnotesize
\begin{tabular}{lll}
\hline
Family & Count & Examples \\
\hline
Complete graphs $K_n$ & 6 & $K_3, K_5, K_7, K_9, K_{11}, K_{13}$ \\
Complete multipartite $K_{n_1,\ldots,n_r}$ & 36 & $K_{1,1,3}, K_{1,3,3}, K_{3,3,3}, K_{1,1,1,1,3}, \ldots$ \\
Strongly regular graphs & 2 & $\text{srg}(9,4,1,2)$, $\text{srg}(10,6,3,4)$ \\
\hline
Total & 44 & \\
\hline
\end{tabular}
\end{table}

The two strongly regular graphs deserve special attention. The graph $\text{srg}(9,4,1,2)$ from $[9,4,4]$ is the Paley graph of order 9, while $\text{srg}(10,6,3,4)$ from $[10,4,4]$ is the complement of the Petersen graph.
According to Grassl's tables~\cite{grassl}, these codes are optimal, as were seen in Examples~\ref{ex:paley} and~\ref{ex:petersen}.

\section{Concluding Remarks and Open Problems}\label{sec:conclusion}

In this paper, we investigated the structural correspondence between LCD codes and simple graphs. We established a complete bijection between binary even LCD codes and simple graphs with idempotent adjacency matrices over $\FF_2$, where permutation equivalence of codes corresponds exactly to graph isomorphism. We characterized which distance-regular graphs yield LCD codes via intersection array parameters. As applications, we provided a theoretical explanation for observation (i) of Haemers et al.\ \cite{haemers1999} on conference graphs and demonstrated that mass formulas for LCD codes can classify graphs with idempotent adjacency matrices.

Several directions merit further investigation. For strongly regular graphs, questions remain about optimality: do Paley graphs always maximize minimum distance among conference graphs? Is the Paley graph the unique optimizer for $\text{srg}(41,20,9,10)$?
More broadly, understanding the relationship between code invariants and graph invariants could reveal deeper connections beyond our idempotency characterization.

\section*{Acknowledgment}

The author thanks Stefka Bouyuklieva for insightful questions and comments that significantly improved this manuscript. This work was supported by JSPS KAKENHI Grant Number JP25K17290.

\end{document}